\documentclass[11pt]{article}
\usepackage{latexsym}
\usepackage{amssymb}
\usepackage{enumerate}
\usepackage{enumitem}
\setlist[enumerate]{leftmargin=0.6cm}
\usepackage{xcolor}
\usepackage{commath}
\usepackage{amsfonts}
\usepackage{amsmath,amsthm}
\usepackage{hyperref}
\usepackage{algorithm}
\usepackage{algorithmic}
\usepackage{framed}
\usepackage{boites}
\usepackage[pdftex]{graphicx} 

\newenvironment{@abssec}[1]{%
    \if@twocolumn

      \section*{#1}%
    \else

      \vspace{.05in}\footnotesize
      \parindent .2in
 {\upshape\bfseries #1. }\ignorespaces
    \fi}

    {\if@twocolumn\else\par\vspace{.1in}\fi}

\newenvironment{keywords}{\begin{@abssec}{\keywordsname}}{\end{@abssec}}
\newenvironment{AMS}{\begin{@abssec}{\AMSname}}{\end{@abssec}}

\newcommand\keywordsname{Key words}
\newcommand\AMSname{AMS subject classifications}
\newcommand\AMname{AMS subject classification}
\newcommand\restr[2]{{
\left.\kern-\nulldelimiterspace 
#1 
\vphantom{|} 
\right|_{#2} 
}}
\newtheorem{theorem}{Theorem}[section]

\newtheorem{remark}[theorem]{Remark}

\newtheorem{problem}{Problem}

\newcommand{\NN}{\mathbb{N}}

\newcommand{\RR}{\mathbb{R}}

\def\XXint#1#2#3{{\setbox0=\hbox{$#1{#2#3}{\int}$}
\vcenter{\hbox{$#2#3$}}\kern-.5\wd0}}

\newcommand{\link}{\mathop{\circ\kern-.35em -}}

\newcommand{\ol}{\overline}
\newcommand{\pa}{\partial}

\newcommand{\dv}{\mathop{\mathrm{div}}}

\newcommand{\gr}{\nabla}

\newcommand{\al}{\alpha}
\newcommand{\be}{\beta}
\newcommand{\ga}{\gamma}  
\newcommand{\Ga}{\Gamma}

\newcommand{\la}{\lambda}

\newcommand{\si}{\sigma}
\newcommand{\Si}{\Sigma}

\newcommand{\Om}{\Omega}

\newcommand{\sg}{\sigma}

\newcommand\setbld[2]{\left\{ #1 \;\middle |\; #2\right\}}

\newcommand{\cC}{\mathcal{C}}

\title{\bf Symmetry breaking solutions for a two-phase overdetermined problem of Serrin-type \thanks{This research was partially supported by the 
Grant-in-Aid for JSPS Fellows No.18J11430 and No.19J12344.}
}

\author{Lorenzo Cavallina ${}^{\dagger}$ and Toshiaki Yachimura \thanks{Research Center for Pure and Applied Mathematics,
Graduate School of Information Sciences, Tohoku
University, Sendai, 980-8579, Japan ({\tt cava@ims.is.tohoku.ac.jp}, {\tt  yachimura@ims.is.tohoku.ac.jp}).}
}
\date{}

\begin{document}
\maketitle
\begin{abstract}{}
In this paper, we consider an overdetermined problem of Serrin-type for a two-phase elliptic operator with piecewise constant coefficients. 
We show the existence of infinitely many branches of nontrivial symmetry breaking solutions
which bifurcate from any radially symmetric configuration satisfying some condition on the coefficients.
\end{abstract}

\pagestyle{plain}
\thispagestyle{plain}

\section{Introduction and main result}\label{introduction}
\begin{figure}[h]
\centering
\includegraphics[width=0.37\linewidth]{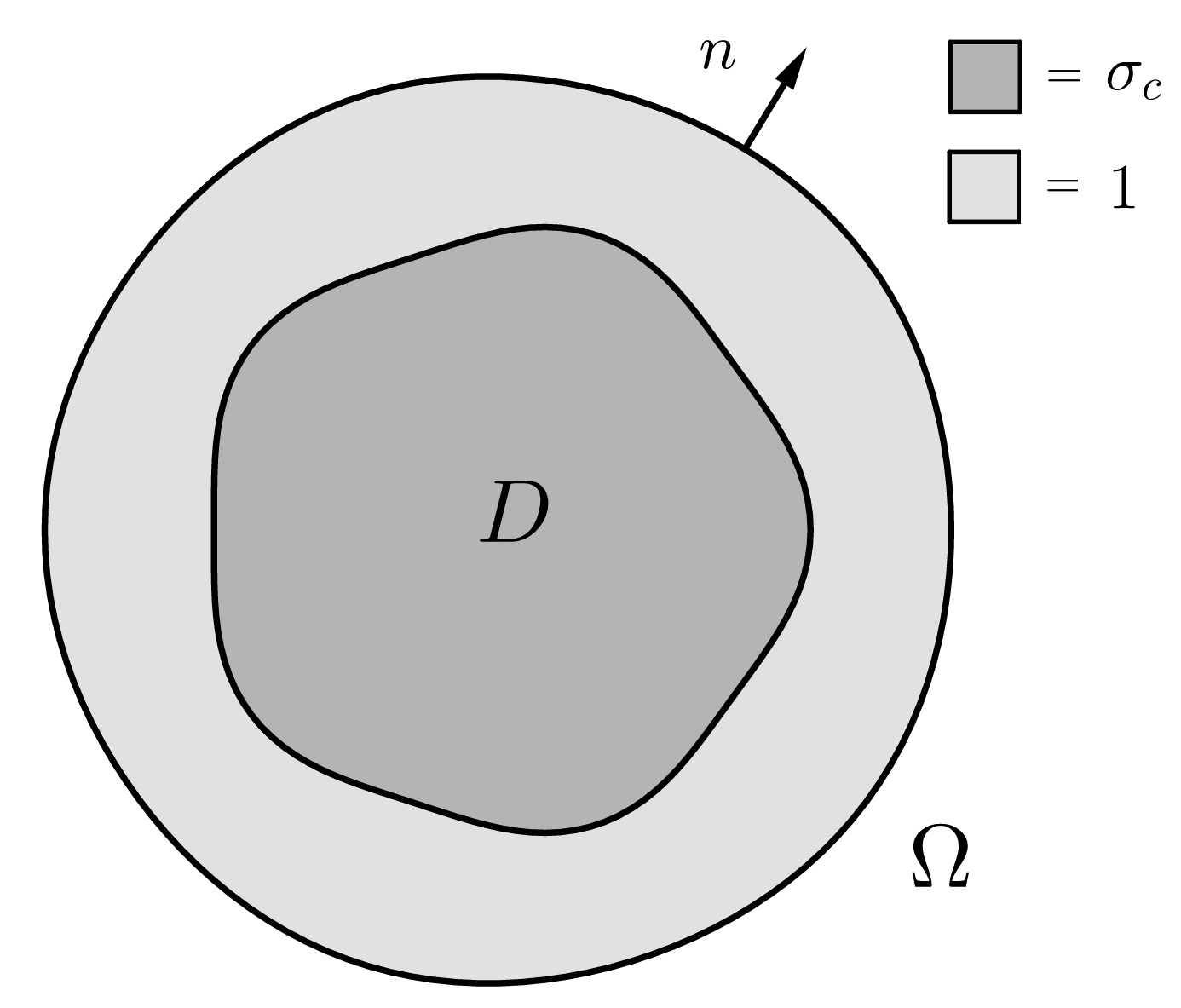}
\caption{Problem setting} 
\label{pb setting}
\end{figure}

In this paper, we consider a bifurcation analysis of a Serrin-type overdetermined problem for an elliptic operator with piecewise constant coefficients. First, let us introduce the problem setting of our overdetermined problem. 
Let $(D,\Omega)$ be a pair of sufficiently smooth bounded domains of $\RR^{N}$ ($N\geq2$) such that $\overline{D} \subset \Omega$. Moreover, let $n$ denote the outward unit normal vector of $\Omega$. 
We consider the following two-phase Serrin-type overdetermined problem:
\begin{equation}\label{odp}
\begin{cases}
-\dv(\sigma \gr u)=1 \quad \textrm{ in }\Omega,\\
u=0\quad \textrm{ on }\partial\Omega,\\
\partial_n u= c \quad \textrm{ on }\partial\Omega, 
\end{cases}
\end{equation}
where $c$ is a real constant and $\sigma = \sigma(D,\Om)$ is the piecewise constant function given by 
\begin{equation*}
\sigma(x) = 
\begin{cases}
\si_c \quad &\text{in} \,\, D, \\
1 \quad &\text{in} \,\, \Omega \setminus D,  
\end{cases}
\end{equation*}
and $\si_c$ is a positive constant such that $\si_c \neq 1$ (Fig. \ref{pb setting}). 

We remark that, if \eqref{odp} is solvable, then the parameter $c$ must be equal to $c(\Omega)=-|\Omega|/|\pa\Omega|$ by integration by parts. 
In what follows, we will say that a pair of domains $(D,\Omega)$ is a solution of problem \eqref{odp} whenever problem \eqref{odp} is solvable for $\sg=\sg(D,\Omega)$. 
Let us define the inner problem and outer problem associated to problem \eqref{odp} (see \cite{CY2019}). 
\begin{problem}[Inner problem]
For a given domain $\Omega$ and a real number $0<V_0<|\Omega|$, find a domain $D\subset\ol D \subset \Omega$ with volume $|D|=V_0$, such that the pair $(D,\Omega)$ is a solution of the overdetermined problem \eqref{odp}.
\end{problem}

\begin{problem}[Outer problem]
For a given domain $D$ and a real number $V_0>|D|$, find a domain $\Omega\supset \ol D$ with volume $|\Omega|=V_0$, such that the pair $(D,\Omega)$ is a solution of the overdetermined problem \eqref{odp}. 
\end{problem}

The case where $D$ is empty (one-phase setting) has been studied by many mathematicians in various situations since the pioneering work of Serrin \cite{Se1971}, who proved that the overdetermined problem \eqref{odp} without the inclusion $D$ is solvable if and only if the domain $\Omega$ is a ball.  
We refer to \cite{BNST2008, BHS2014, MagnaniniPoggesi2017, NT2018} and references therein. 

However, when $D$ is not empty (two-phase setting), there are a few results for the overdetermined problem \eqref{odp}.
The paper \cite{camasa} deals with the inner problem (Problem $1$) of the overdetermined problem \eqref{odp}, the authors proved the local existence and uniqueness for the inner problem near concentric balls. 

The authors, in \cite{CY2019}, proved the following local existence and uniqueness results for the outer problem (Problem $2$) near concentric balls by perturbation arguments by means of shape derivatives and the implicit function theorem for Banach spaces. 
\begin{theorem}\label{mainthm1}
Let us define
\begin{equation}\label{sk}
\begin{aligned}
s(k)&= \frac{k(N+k-1)-(N+k-2)(k-1)R^{2-N-2k}}{k( N+k-1)+k(k-1)R^{2-N-2k}} \text{ for }k = 1,2,\ldots,\\
\Sigma&=\setbld{s\in (0,\infty)}{s=s(k)
\,\text{ for some }k = 1,2,\ldots}.
\end{aligned}
\end{equation}
and let $B_{R} \subset B_{1}$ denote concentric balls of radius $R$ and $1$ respectively. 
If $\si_c \notin \Sigma$, then for every domain $D$ of class $\cC^{2,\alpha}$ sufficiently close to $B_{R}$ in the $\cC^{2,\alpha}$-norm , there exists a domain $\Omega$ of class $\cC^{2,\alpha}$ sufficiently close to $B_{1}$ in the $\cC^{2,\alpha}$-norm such that the outer problem (Problem $2$) admits a solution for the pair $(D,\Omega)$. 
\end{theorem}
\begin{remark}
Notice that, by the definition of $s(k)$ in \eqref{sk}, the quantity $s(k)$ is not necessarily positive for all values of $N$, $k$ and $R$. Indeed, for fixed $N$ and $R$, the quantity $s(k)$ tends to $-1$ as $k\to +\infty$. In particular, this implies that the set $\Sigma$ is finite. 
\end{remark}

From Theorem \ref{mainthm1}, problem \eqref{odp} has a solution near concentric balls except for $\si_c \in \Sigma$. 
Our aim in this paper is to examine the case for $\sigma_c$ near $s(m) \in \Sigma$ in the same situation of Theorem \ref{mainthm1}. In particular, our interest is the shape of the solution of the outer problem near $\si_c \in \Sigma$. 
In what follows, we introduce some notations in order to state the main theorem in this paper precisely. 
\begin{figure}[h]
\centering
\includegraphics[width=0.40\linewidth]{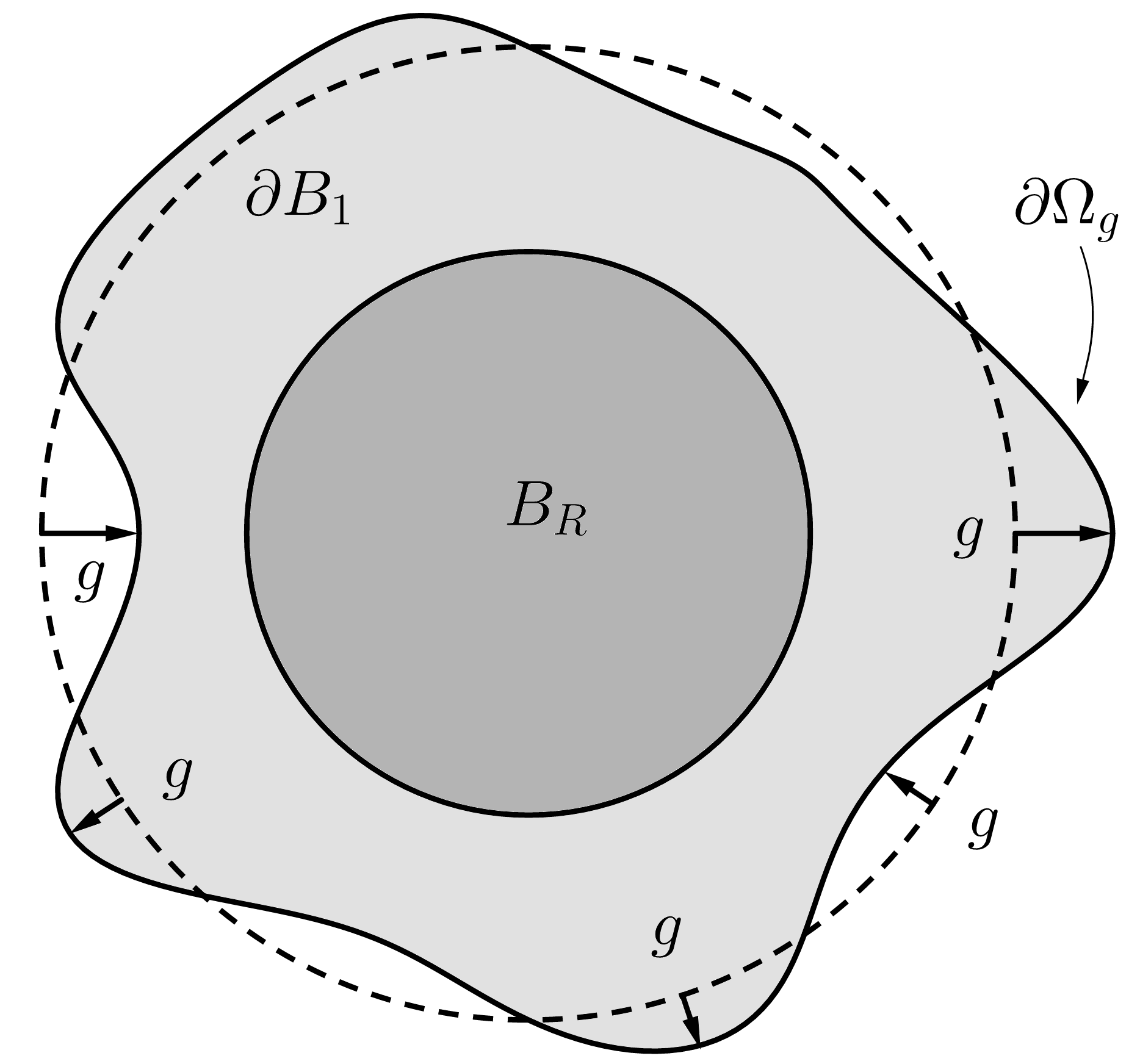}
\caption{The geometrical construction used in the definition of $\Psi(g,\la)$} 
\label{construction}
\end{figure}

Let us take an element $s(m) \in \Si$ for some $m \geq 1$ and let $X$ and $Y$ denote the Banach spaces 
\begin{equation*}
X=\left\{ g\in \cC^{2,\al}(\pa B_1) \;:\; \int_{\pa B_1} g=0\right\}, \quad
Y=\left\{ h\in \cC^{1,\al}(\pa B_1) \;:\; \int_{\pa B_1} h=0\right\}
\end{equation*}
endowed with their natural norms. We consider the functional $\Psi: X\times \RR\to Y$ defined by
\begin{equation}\label{Psi}
    \Psi(g,\la) = \left\lbrace \pa_{n_g} v_{g} - c_g  \right\rbrace J_\tau (g). 
\end{equation}
In what follows, we will explain the notation involved in the definition of \eqref{Psi}. 
For $g \in X$, let $\Om_g$ be the unique bounded domain whose boundary is defined as
\begin{equation*}
    \pa\Om_g=\left\{ x+g(x)n(x) \;:\; x\in\pa B_1\right\}
\end{equation*}
with outward unit normal vector denoted by $n_g$.
Moreover, let $v_g$ be the solution of the Dirichlet boundary value problem given by the first two equations in \eqref{odp} for $(D,\Om)=(B_R,\Om_g)$ and $\sg_c=s(m)+\la$. By the definition \eqref{Psi}, we notice that $\Psi(0,\la) = 0$ for any $\la$ since the pair of concentric balls $(B_R,\Omega_0)$ is a solution of overdetermined problem \eqref{odp}. 
By a slight abuse of notation, we will use $\pa_{n_g}v_g$ to denote the function of value 
\begin{equation*}
    \pa_{n_g}v_g\left( x+g(x)n(x)\right) \quad \text{for }x\in\pa B_1.
\end{equation*}
Finally, $c_g=c(\Om_g)=-|\Om_g|/|\pa\Om_g|$ and $J_\tau(g)$ denotes the tangential Jacobian of the map $x\mapsto x+g(x)n(x)$ from $\pa B_1$ to $\pa\Om_g$ (see \cite[Definition 5.4.2]{HP2005}). 

Now we can present the main result in this paper. 
\begin{theorem}\label{mainthm2}
Let $(D,\Om)=(B_R,\Om_g)$ $(0 < R < 1)$. Also we take an element $s(m) \in \Si$ for some $m \geq 1$ and suppose that $\si_c = s(m) + \la$, where $\la \in \RR$. If we consider the equation 
\begin{equation*}
    \Psi(g,\la) = 0,  
\end{equation*}
then $(0,0)$ is a bifurcation point of the equation $\Psi(g,\la) = 0$. That is, there exists a smooth function $\varepsilon\mapsto \la(\varepsilon)\in\RR$ with $\la(0)=0$ 
such that overdetermined problem \eqref{odp} admits a nontrivial solution of the form $(B_R,\Om_{g(\varepsilon)})$ for $\sg_c=s(m)+\la(\varepsilon)$ and $\varepsilon$ small (notice that, by continuity, $\sg_c=s(m)+\la(\varepsilon)>0$ if $\varepsilon$ is small enough). If $N=2$, then the symmetry breaking solution $(B_R,\Om_{g(\varepsilon)})$ satisfies 
\begin{equation}\label{symmetry breaking solutions}
    g(\varepsilon) = \varepsilon \cos(m\theta) + o(\varepsilon) \quad \text{in }\cC^{2,\al}(\pa\Om_0)\quad \text{as }\varepsilon\to0.   
\end{equation}
Moreover, if $N \ge 3$, then there exists a spherical harmonic $Y_m$ of $m$-th degree, such that the symmetry breaking solution $(B_R,\Om_{g(\varepsilon)})$ satisfies 
\begin{equation}\label{symmetry breaking solutions 2}
    g(\varepsilon) = \varepsilon Y_m(\theta) + o(\varepsilon) \quad \text{in }\cC^{2,\al}(\pa\Om_0) \quad \text{as }\varepsilon\to0.   
\end{equation}

\end{theorem}
From Theorem \ref{mainthm2}, if $D = B_R$, then the outer problem has solutions not only for $\Omega = B_1$ but also for $\Omega=\Omega_g$ given by \eqref{symmetry breaking solutions} and \eqref{symmetry breaking solutions 2}. That is, there exist branches of symmetry breaking solutions of the outer problem emanating from the bifurcation points $\si_c \in \Sigma$. This implies that the uniqueness of the outer problem does not hold near $\si_c \in \Sigma$ because symmetry breaking phenomena occur. 
Similar results appear in the context of free boundary problems of a circulating flow with surface tension \cite{Oka1984} and a model of tumor growth \cite{FR2001, EM2011}.  

This paper is organized as follows. In Section \ref{bifur}, we prove Theorem \ref{mainthm2} when $N = 2$. This proof is based on the results obtained in \cite{CY2019} and the Crandall--Rabinowitz theorem. In Section \ref{bifur2}, we consider high dimensional case $N \ge 3$ and establish Theorem \ref{mainthm2}. 

\section{Proof of Theorem \ref{mainthm2} for $N = 2$}\label{bifur}
In this section, we prove the main theorem \ref{mainthm2}. We obtain the existence of symmetry breaking bifurcation solutions of  overdetermined problem \eqref{odp} applying the following version of the Crandall--Rabinowitz theorem. 
\begin{theorem}[Crandall--Rabinowitz theorem \cite{CR1971}]\label{Crandall--Rabinowitz theorem}
Let $X$, $Y$ be real Banach spaces and $\Psi(x,\la)$ be a $C^p$ map ($p \geq 3$) of a neighborhood $(0,\la_0)$ in $X \times \RR$ into $Y$. Suppose that 
\begin{enumerate}[label=(\roman*)]
\item $\Psi(0,\la) = 0$ for all $\la$ in a neighborhood of $\la_0$. 
\item There exists $x_0 \in X$ such that {\rm Ker}\,$\pa_x\Psi(0,\la_0)$ is a one-dimensional space spanned by $x_0$. 
\item {\rm Im}\,$\pa_x\Psi(0,\la_0)$ is a closed subspace of $Y$ which has codimension one. 
\item $\pa_\la\pa_x\Psi(0,\la_0)[x_0] \notin$ {\rm Im}\,$\pa_x\Psi(0,\la_0)$.
\end{enumerate}
Then $(0, \la_0)$ is a bifurcation point of the equation $\Psi(x,\la)=0$ in the following sense: In a neighborhood of $(0, \la_0)$ the set of solutions of $\Psi(x,\la) = 0$ consists of two $C^{p-2}$ smooth curves $\Gamma_1$ and $\Gamma_2$ which intersect only at the point $(0,\la_0)$; $\Gamma_1$ is the curve $(0,\la)$ and $\Gamma_2$ can be parametrized as follows: 
\begin{equation*}
\Gamma_2: \left( x(\varepsilon),\la(\varepsilon) \right), \quad \varepsilon:\text{small}, \quad \left( x(0),\la(0) \right) = (0,\la_0), \quad x'(0)=x_0.  
\end{equation*}
\end{theorem}
In what follows, we assume that $N=2$. 
\begin{proof}[Theorem \ref{mainthm2}, $N=2$]
Take an element $s(m)\in\Si$ for some $m \geq 1$ and let $\Psi$ be the map defined by \eqref{Psi}. By definition, notice that $\Psi(g,\la)=0$ if and only if the pair $(B_R,\Om_g)$ solves \eqref{odp} for $\sg_c=s(m)+\la$.
By \cite[Theorem 3.15 (iii)]{cavaphd}, the map $\Psi$ is Fr\'echet differentiable infinitely many times in a neighborhood of the origin in $X$. Moreover, by the explicit formula of its Fr\'echet derivative $\pa_x \Psi(0,\la)$ computed in \cite[Theorem 3.5]{CY2019}, we know that ${\rm Ker}\ \pa_x\Psi(0,\la)$ is a two dimensional space, spanned by $\{\cos(m\theta),\sin(m\theta)\}$. As a consequence, we cannot apply the Crandall--Rabinowitz theorem (Theorem \ref{Crandall--Rabinowitz theorem}) directly. In order to reduce the kernel to a one dimensional space, we introduce the following spaces of even functions:
\begin{align*}
    X^* &=\left\{ g\in X \;:\; g(\theta)=g(2\pi-\theta),\ \theta\in [0,2\pi) \right\},\\
Y^* &=\left\{ h\in Y \;:\; h(\theta)=h(2\pi-\theta) ,\ \theta\in [0,2\pi) \right\},
\end{align*}
where we identified the unit circle $\pa B_1\subset \RR^2$ with the interval $[0,2\pi)$.
Now, we consider the restriction $\Psi^*$ of $\Psi$ on $X^*$. We claim that $\Psi^*$ is a well-defined mapping 
\begin{equation*}
    \Psi^*: X^*\to Y^*.
\end{equation*}
To show this, notice that $g\in X^*$ implies that the configuration $(B_R, \Om_g)$ is symmetric with respect to the $x$-axis. Now, by the unique solvability of the Dirichlet boundary value problem given by the first two equations in \eqref{odp}, this implies that also $v_g$ shares the same symmetry and, thus, $\Psi^*(g,\la)=\Psi(g,\la)\in Y^*$ as claimed.

We will now apply Theorem \ref{Crandall--Rabinowitz theorem} to the map $\Psi^*$.
Recall that $\Psi^*(0,\la) = 0$ for any $\la$ since the pair of concentric balls $(B_R,\Omega_0)$ is a solution of overdetermined problem \eqref{odp}. 
This fact implies that $(i)$ holds true. 
Let us check condition $(ii)$. In the proof of Theorem $3.6$ in \cite{CY2019}, we computed the Fr\'echet derivative $\pa_x\Psi(0,\la)$. The case $N=2$ reads
\begin{equation}\label{preserves eigenspaces}
\pa_x\Psi(0,\la)[g]= \sum_{k=1}^\infty \beta_k(\la)\left(\alpha_{k}^{\rm even}
\cos(k\theta)+\al_k^{\rm odd}\sin(k\theta)       \right), 
\end{equation}
for
\begin{equation*}
g=\sum_{k=1}^\infty 
\left(
\al_k^{\rm even} \cos(k\theta) 
+\al_k^{\rm odd}\sin(k\theta)
\right), 
\end{equation*}
where
\begin{equation}\label{beta_k}
\beta_k(\la)= \frac{(k+1)(s(m)+\la-1)k+(k+ks(m)+k\la)(k-1)R^{-2k}}{2(k+ks(m)+k\la)R^{-2k}+2k  (1-s(m)-\la)}.
\end{equation}
Now, a simple computation with \eqref{sk} at hand yields that 
\begin{equation*}
    \be_m(0)=0 \quad \text{and}\quad \be_k(0)\ne 0 \ \text{for }k\ne m.
\end{equation*}
Let $x_0 = cos(m\theta)$. 
Notice that $X^*$ is the subspace of $X$ spanned by $\left\{\cos(k\theta)\right\}_{k\ge1}$. Then, combining \eqref{preserves eigenspaces} with the fact that $\be_m(0)=0$ by construction, we obtain
\begin{equation*}
    \text{Ker}\ \pa_x\Psi^*(0,0) = \text{span} \{ x_0 \}. 
\end{equation*}
Thus condition $(ii)$ holds true. 
Moreover, $\pa_x\Psi^*(0,0)[\cos(k\theta)] = \beta_k(0) \cos(k\theta)$, where $\beta_k(0) \neq 0$ for $k \neq m$. This implies that \begin{equation*}
    \text{Im}\ \pa_x\Psi^*(0,0) \oplus \text{Ker}\ \pa_x \Psi^*(0,0) = Y^*. 
\end{equation*}
Therefore, $\text{Im}\ \pa_x\Psi^*(0,0)$ is codimension one and thus also condition $(iii)$ holds true. 

Let us finally check condition $(iv)$. Note that 
\begin{equation*}
    \pa_x\Psi^*(0,\la)[x_0] = \beta_m(\la)\  x_0.  
\end{equation*}
By \eqref{beta_k}, we can easily compute that 
\begin{equation}\label{dif beta}
\partial_{\la} \beta_m (0) = \dfrac{\{ m(m+1) + m(m-1)R^{-2m}\} \{ 2(m+ms(m))R^{-2m}+2m(1-s(m))\}}{\{2m(1+s(m))R^{-2m} + 2m(1-s(m))\}^2}, 
\end{equation}
where we used the fact that $\beta_m(0) = 0$ by construction. Since $0 < s(m) < 1$ and $m \geq 1$, the right hand side of \eqref{dif beta} is positive. Thus by \eqref{preserves eigenspaces}, 
\begin{equation*}
    \pa_\la\pa_x\Psi^*(0,0) [x_0] = \partial_{\la} \beta_m (0) \, x_0 \in \text{Ker}\  \pa_x\Psi^*(0,0)\setminus\{0\}. 
\end{equation*}
In other words, 
\begin{equation*}
\pa_\la\pa_x\Psi^*(0,0) [x_0] \notin \text{Im}\ \pa_x\Psi^*(0,0).  
\end{equation*}

Therefore, by the Crandall--Rabinowitz theorem (Theorem \ref{Crandall--Rabinowitz theorem}), $(0,0)$ is a bifurcation point of the equation $\Psi(g,\la) = 0$ in the sense that there exists a $\cC^\infty$ curve $(g(\cdot),\la(\cdot))$ from a neighborhood of $0\in\RR$ into $X^*\times \RR$, with $(g(0),\la(0))=(0,0)$ and such that, for all small $\varepsilon$, there exists a symmetry breaking solution of overdetermined problem \eqref{odp} for $\sg_c=s(m)+\la(\varepsilon)$, represented by $(B_R,\Om_{g(\varepsilon)})$, with 
\begin{equation*}
    g(\varepsilon) = \varepsilon \cos(m\theta) + o(\varepsilon) \quad \text{as }\varepsilon\to0.   
\end{equation*}
\end{proof}
\begin{remark}
Theorem \ref{mainthm2} ensures the existence of nontrivial solutions of \eqref{odp} of the form $(B_R,\Om)$. In particular, such solutions only partially inherit the symmetry of the core $B_R$. 
One might wonder whether nontrivial solutions of the form $(D,B_1)$ exist for some subdomain $D$ other than a ball. Actually this is not the case, since Theorem $5.1$ of \cite{Sakaguchi Bessatsu} states that, if $B_1\setminus D$ is connected and the pair $(D,B_1)$ solves \eqref{odp}, then $D$ and $B_1$ must be concentric balls.
\end{remark}
\begin{figure}[h]
\centering
\includegraphics[width=0.7\linewidth]{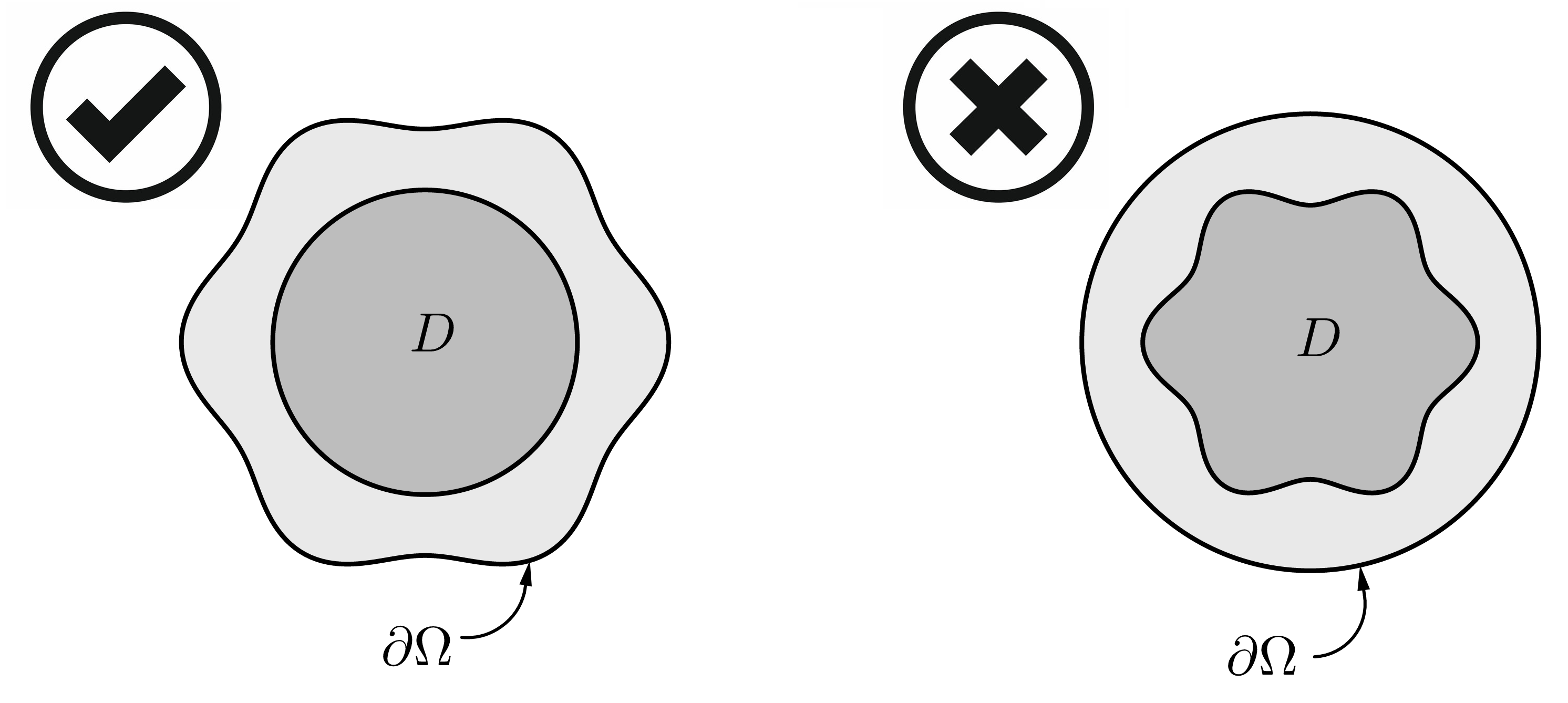}
\caption{Left: a symmetry-breaking bifurcating solution of \eqref{odp} given by Theorem \ref{mainthm2} ($m=6$).\newline Right: a symmetry-breaking configuration that cannot be a solution to \eqref{odp} in light of Theorem 5.1 of \cite{Sakaguchi Bessatsu}} 
\label{nonlinear}
\end{figure}

\section{Proof of Theorem \ref{mainthm2} for $N\ge3$}\label{bifur2}

The proof of Theorem \ref{mainthm2} when $N\ge3$ follows along the same lines as the previous section. Indeed, as in the case $N=2$, the Crandall--Rabinowitz theorem cannot be applied directly because ${\rm{Ker}}\, \pa_x\Psi(0,0)$ is not one dimensional.
By the $N$-dimensional analogous of \eqref{preserves eigenspaces} (see \cite[equation (3.12)]{CY2019}), ${{
\rm Ker}\, \pa_x \Psi(0,0)}$ is the subspace of $X$ spanned by the spherical harmonics whose degree is $m$.
In order to reduce the kernel to a one dimensional space, we follow the same ideas as \cite{KS2019} and consider the restriction $\Psi^*$ of $\Psi$ to the space $X^*$ of functions in $X$ that are invariant with respect to some specific group of symmetries $\Ga\subset O(N)$.
Here we recall the definition of $\Ga$-invariance with respect to a subgroup $\Ga$ of the orthogonal group $O(N)$. 
A function $g\in X$ is said to be $\Ga$-invariant if 
\begin{equation*}
    g(\theta)=g(\ga(\theta)) \quad \text{for all }\theta\in\pa B_1,\ \ga\in\Ga.
\end{equation*}
If, for example, we set $\Ga={\rm Id}\times O(N-1)$, then the space of $\Ga$-invariant spherical harmonics of any given degree $k\in \NN$ is a one dimensional space. 
In particular, if $X^*\subset X$ is the subset of $\Ga$-invariant functions and $\Psi^*$ the restriction of $\Psi$ to $X^*$, then also ${\rm Ker} \, \pa_x \Psi^*(0,0)$ is a one dimensional space, which can be considered to be spanned by some spherical harmonic $x_0\in X^*$. The rest of the proof runs just as the one in Section \ref{bifur}, by checking conditions $(i)-(iv)$ of the Crandall--Rabinowitz theorem applied to the map $\Psi^*$.

\begin{remark}
We claim that all nontrivial solutions $(B_R,\Om_{g(\varepsilon)})$ given by Theorem \ref{mainthm2} share the same symmetries of the element $x_0\in X^*$, defined such that ${\rm Ker} \, \pa_x \Psi^*(0,0)={\rm span}\{x_0\}$. To this end, let $\Ga\subset O(N)$ be a symmetry group such that the function $x_0$ is $\Ga$-invariant. Now, consider the further restriction $\Psi^{**}$ of $\Psi^*$ to the subspace $X^{**}$ of all $\Ga$-invariant functions in $X^*$. Notice that, since $x_0$ is $\Ga$-invariant by hypothesis, then ${\rm Ker} \, \pa_x\Psi^*(0,0)={\rm Ker} \, \pa_x \Psi^{**}(0,0)={\rm span}\{x_0\}$. Another application of the Crandall--Rabinowitz theorem to $\Psi^{**}$ yields that $g(\varepsilon)$ is also $\Ga$-invariant. The claim follows by the arbitrariness of $\Ga$.
\end{remark}

\section*{Acknowledgements}
We would like to sincerely thank the anonymous reviewers for their comments and observations, which helped to improve the overall readability of the paper.

\begin{small}

\end{small}
\end{document}